\documentclass[12pt,reqno]{amsart}

\usepackage[final
]{showkeys}

\usepackage{AKstyle}

\usepackage[table]{xcolor}

\newcommand{\vertiii}[1]{{\left\vert\kern-0.25ex\left\vert\kern-0.25ex\left\vert #1
    \right\vert\kern-0.25ex\right\vert\kern-0.25ex\right\vert}}

\usepackage{tikz}
\usepackage{enumerate}

\usepackage{caption}
\usepackage[labelfont=rm]{subcaption}

\usepackage[pagebackref=true, colorlinks]{hyperref}

\hypersetup{pdffitwindow=true,linkcolor=blue,citecolor=blue,urlcolor=blue,menucolor=blue}

\usepackage{comment}

\newcommand{\StructThm}{\hyperref[thm:tildeC]{Structure Theorem}}

\usepackage{arydshln}

\begin{document}

\author{Alex Kontorovich}
\thanks{Kontorovich is partially supported by
an NSF CAREER grant DMS-1455705, an NSF FRG grant DMS-1463940,  a BSF grant, a Simons Fellowship,  a von Neumann Fellowship at IAS, and the IAS's NSF grant DMS-1638352.}
\email{alex.kontorovich@rutgers.edu}
\address{Department of Mathematics, Rutgers University, New Brunswick, NJ, and School of Mathematics, Institute for Advanced Study, Princeton, NJ}
\author{Kei Nakamura}
\thanks{Nakamura is partially supported by
 NSF FRG grant DMS-1463940.}
\email{kei.nakamura@rutgers.edu}
\address{Department of Mathematics, Rutgers University, New Brunswick, NJ}

\title
[Geometry 
and Arithmetic of
Crystallographic
Sphere Packings]
{
Geometry 
and Arithmetic of\\ 
Crystallographic
Sphere Packings}

\begin{abstract}
We introduce the notion of a ``crystallographic sphere packing,'' defined to be one whose limit set is that of a geometrically finite hyperbolic reflection group in one higher dimension.
We  exhibit for the first time an infinite family of conformally-inequivalent such
with all radii being reciprocals of integers.
We then prove a result in the opposite direction:  the ``superintegral'' ones exist only in finitely many ``commensurability classes,'' all in dimensions below $30$.
\end{abstract}
\date{\today}
\maketitle



The
goal  of this program, the details of which will appear elsewhere, 
is to understand the basic ``nature'' of the classical Apollonian gasket. 
Why does its integral structure exist
? (Of course it follows 
here
from Descartes' Kissing Circles Theorem, but is there  
a more fundamental, intrinsic explanation?) Are there more like it? (Around a half-dozen similarly integral circle and sphere packings were previously known, each given by an ad hoc description.)  If so, how many more? Can they be classified?
We develop a basic unified framework for addressing these questions, and find two surprising 
(and 
opposing)
phenomena:
\begin{enumerate}[(I)]
\item\label{I}  there is indeed a whole infinite zoo of integral sphere packings, and 
\item\label{II} up to ``
commensurability,'' there are only {\it finitely-many} Apollonian-like objects, over all dimensions! 
\end{enumerate}


\begin{Def}
By an $S^{n-1}$-{\it packing} (or just  {packing}) $\sP$ 
of 
 $\widehat{\R^{n}}:=\R^{n}\cup\{\infty\}$, we mean an infinite collection of {oriented} $(n-1)$-spheres (or  co-dim-1 planes)  so that:
 \begin{itemize}
 \item the interiors of spheres are disjoint,
 and
 \item the spheres
  densely fill up space; that is, we require that
any ball in $\widehat{\R^{n}}$ intersects the interior of some sphere in $\sP$. 
\end{itemize}
The 
 {\it bend} of a sphere is the reciprocal of its (signed) radius.\footnote{In dimensions $n=2$, that is, for circle packings, the bend is just the curvature. But in higher dimensions $n\ge3$, the various ``curvatures''  
of an $(n-1)$-sphere are proportional to $1/$radius$^{2}$, not $1/$radius; so we instead use the term ``bend''.}
To be dense but disjoint, the spheres in the packing $\sP$ must have arbitrarily small radii, so arbitrarily large bends.
If every sphere in $\sP$ has integer bend,
then we call the  packing  {\it integral}. 
\end{Def}

Without more structure, one can make completely arbitrary constructions of integral packings. 
A key property enjoyed by the classical Apollonian circle packing and connecting it to the theory of ``thin groups'' (see \cite{Sarnak2014, Kontorovich2014}) is that
 it arises as
  the limit set of a geometrically finite reflection group in  hyperbolic space of one higher dimension.

\begin{Def}
We call a 
packing $\sP$
   {\it 
   crystallographic}
if 
its limit set is that of
some
 geometrically finite
 reflection group
 $\G<\text{Isom}(\bH^{n+1})$.
\end{Def}

This definition is sufficiently general to encompass all previously proposed generalizations of Apollonian gaskets found in the literature, including \cite{Boyd1974, Maxwell1982, ChenLabbe2015, GuettlerMallows2010, ButlerGrahamGuettlerMallows2010, Stange2015, Baragar2017}.
With these two basic and general definitions in place, we may already state our first main result, confirming \eqref{I}. 

\begin{thm}\label{thm:infMany}
There exist infinitely many conformally-inequivalent integral crystallographic packings.
\end{thm}

\begin{figure}
\includegraphics[width=.6\textwidth]{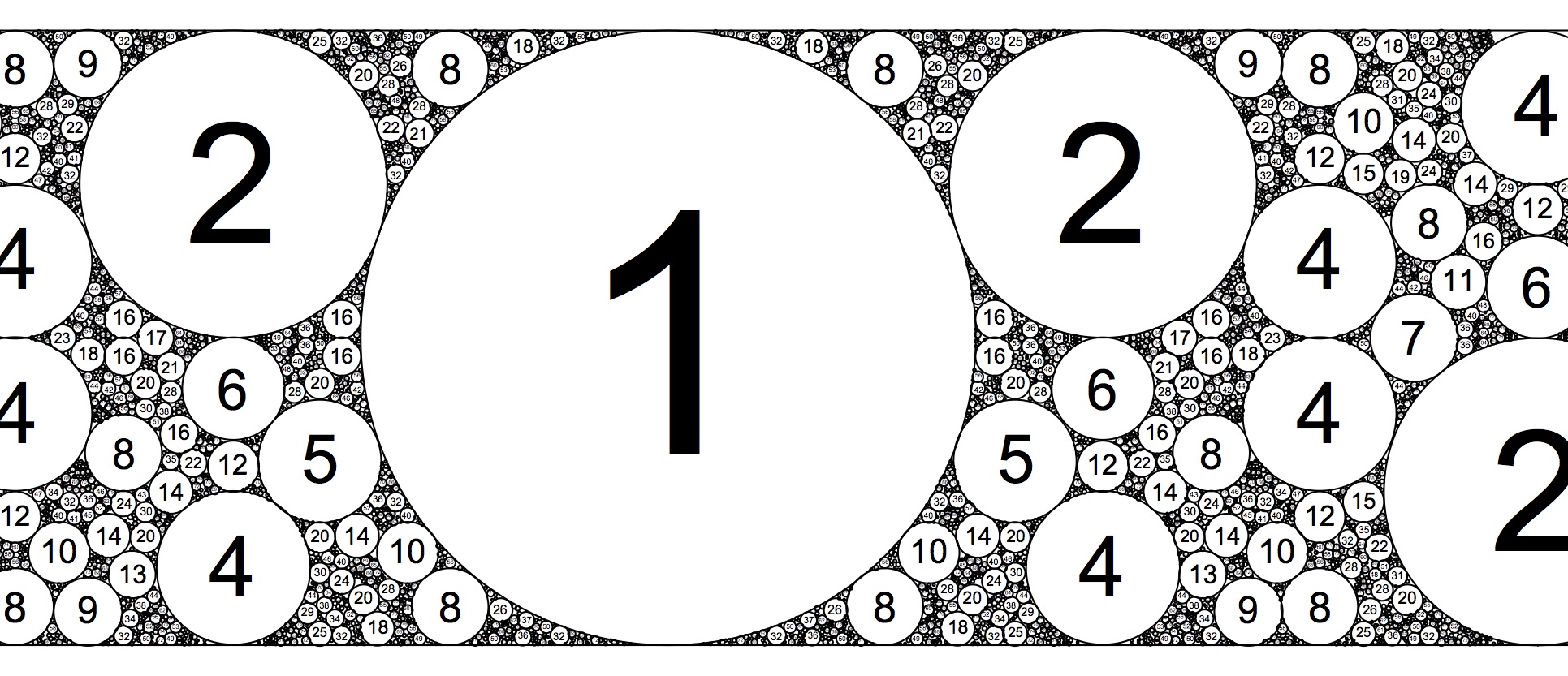}
\caption{A new integral crystallographic packing.  (The circles are labeled with their bend.)}
\label{fig:newPack}
\end{figure}

We show in \figref{fig:newPack} but one illustrative new example,  whose only ``obvious'' symmetry is 
a central mirror image.
It turns out (but may be hard to tell just from the picture) that this packing does indeed arise as the limit set of a Kleinian reflection group.
The argument leading to \thmref{thm:infMany} comes from constructing circle packings ``modeled on'' combinatorial types of convex polyhedra, as follows.
\\

\noindent
{\bf \S(\ref{I}): Polyhedral Packings}\label{sec:1.1}\

\begin{figure}
        \begin{subfigure}[t]{0.25\textwidth}
                \centering
		\includegraphics[width=\textwidth]{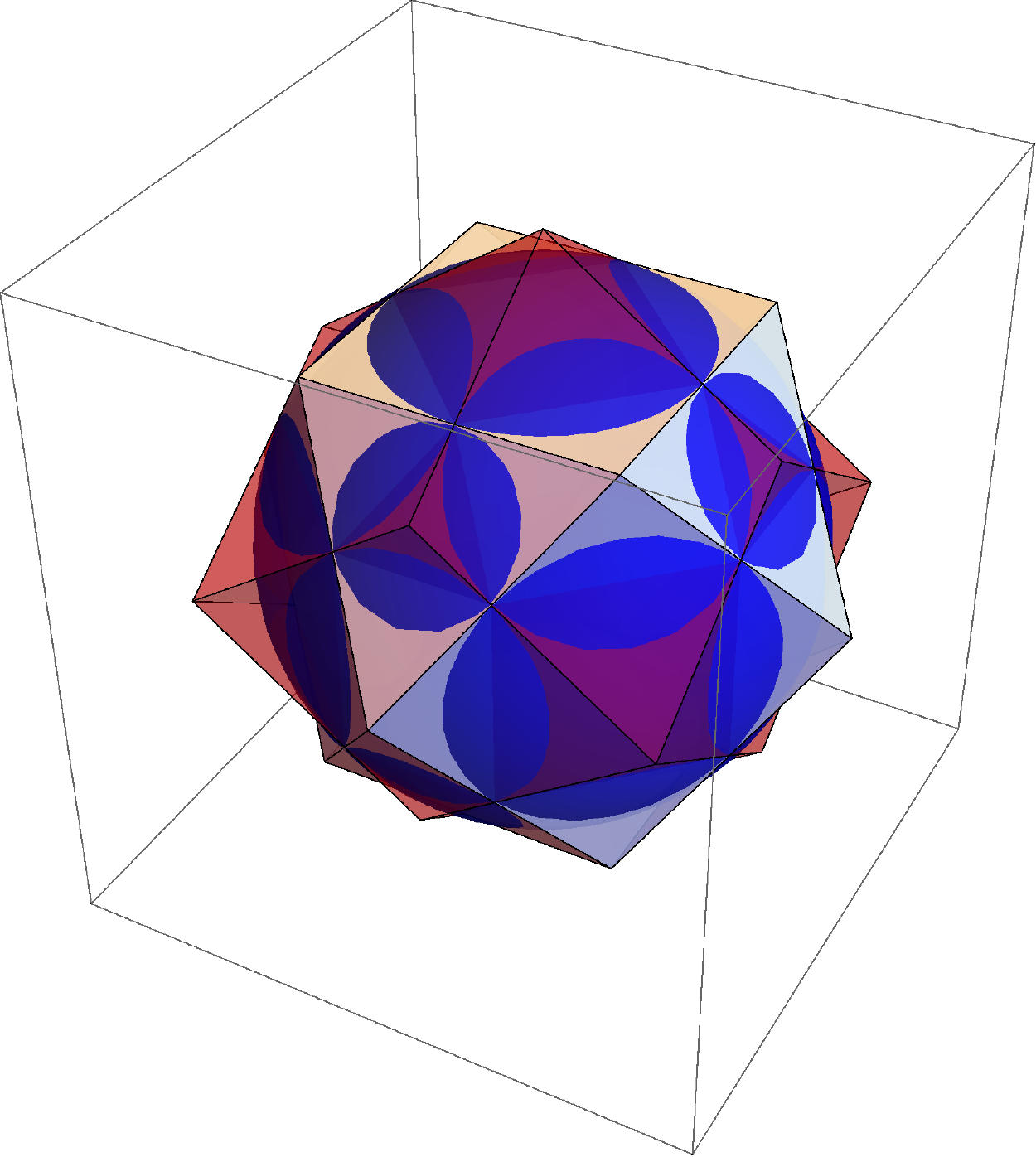}
                \caption{$\Pi$, $\widehat\Pi$ and their midsphere}
                \label{fig:midSph}
        \end{subfigure}%
\qquad
        \begin{subfigure}[t]{0.65\textwidth}
                \centering
		\includegraphics[width=\textwidth]{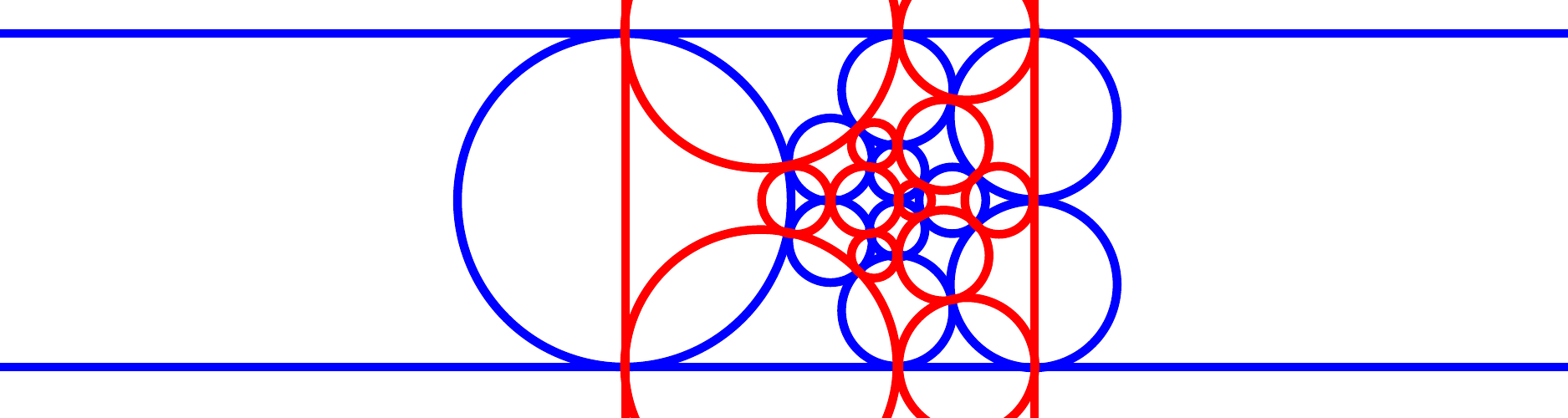}
\caption{An orthogonal cluster and cocluster pair with nerves $\Pi$ and  $\hat\Pi$} 
\label{fig:geom}
        \end{subfigure}

\caption{Geometrization, and cluster/cocluster pair for $\Pi=$ cuboctahedron with dual $\widehat\Pi=$ rhombic dodecahedron.}
\label{fig:}
\end{figure}

Let $\Pi$ be a combinatorial type of a convex polyhedron. 
Equivalently, $\Pi$ is a 3-connected\footnote{Recall that a graph is $k$-connected if it remains connected whenever fewer than $k$ vertices are removed.} planar graph.
A version of the Koebe-Andreev-Thurston Theorem says that there exists a geometrization of $\Pi$ (that is, a realization of its vertices in $\R^3$ with straight lines as edges and faces contained in Euclidean planes) having a midsphere (meaning, a sphere tangent to all edges). 
This midsphere is then also simultaneously a midsphere for the dual polyhedron $\widehat\Pi$.
\figref{fig:midSph} shows the case of a cuboctahedron and its dual, the rhombic dodecahedron.


Stereographically projecting to $\widehat{\R^2}$, we obtain a {\it cluster} (just meaning, a finite collection) $\cC$ of circles whose nerve (that is, tangency graph) is isomorphic to $\Pi$, and a {\it cocluster}, $\widehat\cC$, with nerve $\widehat\Pi$ which meets $\cC$ orthogonally. Again, the example of the cuboctahedron is shown in \figref{fig:geom}.


\begin{Def}
The orbit 
$\sP=\sP(\Pi)=\G\cdot\cC$ 
of the cluster $\cC$ under the group $\G=\<\widehat\cC\>$ generated by reflections through the cocluster $\widehat\cC$
 is
  said to be {\it modeled on} the polyhedron $\Pi$. 
\end{Def}

\begin{lem}
An orbit modeled on a polyhedron is 
a crystallographic 
packing. 
\end{lem}

\begin{figure}
\includegraphics[width=.7\textwidth]{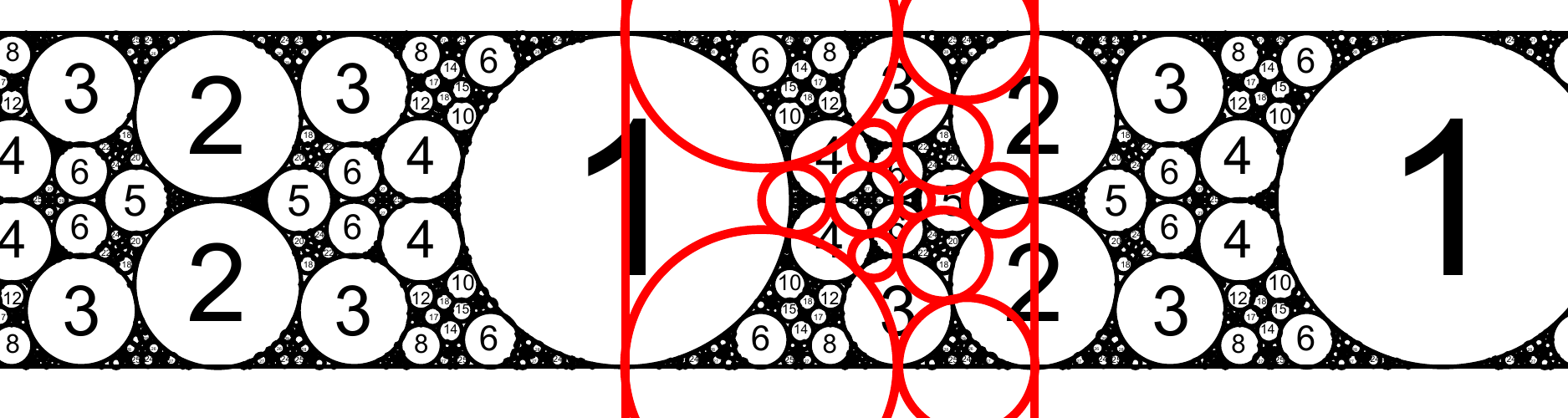}
\caption{A packing modeled on the cuboctahedron, shown with cocluster}
\label{fig:cuboctPack}
\end{figure}

See \figref{fig:cuboctPack} for a packing modeled on the cuboctahedron.
Such packings are unique up to conformal/anticonformal maps by Mostow rigidity, but M\"obius transformations do not generally preserve arithmetic.

\begin{Def}
We call a polyhedron $\Pi$ {\it integral} if there exists {\it some} packing 
modeled on $\Pi$ which is integral. 
\end{Def}

It is not hard to see that the cuboctahedron is indeed integral. So is the tetrahedron, which corresponds to the classical Apollonian gasket.
It is a  fundamental problem to classify all integral polyhedra.

Let us point out some basic  difficulties with this problem. First of all, it is non-trivial to determine whether, given a particular polyhedron, there exists some packing modeled on it which is integral. Indeed, Koebe-Andreev-Thurston geometrization is an infinite limiting process, and how is one to know whether 3.9999 is really 4? To the rescue is Mostow rigidity again, which implies that one can always find cluster/cocluster configurations with all centers and radii {\it algebraic}. 
This means that after computing enough decimal places, one can {\it guess} what the nearby algebraic values might be, and then {\it rigorously} verify whether the guess gives the correct tangency data. This algorithm works for small examples, but once $\Pi$ is sufficiently complicated, it may take a very long time for the guessing process to halt. 
\\

Despite these difficulties, we are able to show the following towards (\ref{I}).

\begin{thm}\label{thm:infInt}
Infinitely many polyhedra are integral, and give rise to infinitely many conformally-inequivalent integral polyhedral packings.
\end{thm}

This of course implies  \thmref{thm:infMany}. To explain the main ideas in the proof, we need some more notation.
\\

Returning to the general setting of crystallographic packings, 
recall that $\sP$ is assumed to arise as the limit set of a discrete group $\G$;
we call the latter
 a {\it symmetry group} of $\sP$.

\begin{Def}
Given a 
packing $\sP$ with 
symmetry group $\G$, we define its {\it supergroup}, $\widetilde\G$, 
to be
 the group generated by $\G$ itself, plus  reflections 
 through all spheres in $\sP$. 
Abusing notation, we may write this as
$$
\widetilde\G \ : = \ \<\G,\sP\> \ < \ \text{Isom}(\bH^{n+1}).
$$
\end{Def}

In the case of a polyhedral packing $\sP=\sP(\Pi)$, the supergroup is simply the group generated by reflections in both the cluster and cocluster, $\widetilde\G=\<\widehat\cC,\cC\>$.

\begin{Def}
The {\it superpacking}, $\widetilde \sP$, of $\sP$ with symmetry group $\G$
is the orbit of $\sP$ under its supergroup, that is,
$$
\widetilde \sP \ :  = \ \widetilde\G\cdot\sP.
$$
\end{Def}

Note that the superpacking is {\it not} a packing by our definition as the sphere interiors are no longer disjoint.\footnote{A related notion of superpacking for the classical Apollonian gasket arose already in work of Graham-Lagarias-Mallows-Wilks-Yan \cite{GLMWYII}; see also the viewpoint of ``Schmidt arrangements'' in work of Stange \cite{Stange2015} and Sheydvasser \cite{Sheydvasser2017}.}

\begin{Def}
We call a 
packing $\sP$ {\it superintegral} if every bend in its superpacking $\widetilde\sP$ is 
an integer.\footnote{Note that 
an unrelated notion of ``superintegrality'' is defined
in \cite[\S8]{GLMWYII}.}
\end{Def}

\begin{rmk}
While different symmetry groups $\G$ lead to different (but commensurate) supergroups $\widetilde\G$, the superpackings are universal, the same for all choices of $\G$.
\end{rmk}

Returning to polyhedral packings, we say that a polyhedron is {\it  superintegral} if some packing modeled on it is.
To prove \thmref{thm:infInt}, we actually prove the following stronger statement.

\begin{thm}\label{thm:infSuper}
Infinitely many polyhedra are superintegral, and give rise to infinitely many conformally-inequivalent superintegral crystallographic packings.
\end{thm}

\begin{rmk}
We stress the conformal-inequivalence here because it turns out that infinitely many polyhedra give rise to the {\it same} crystallographic packing; so the
first part of \thmref{thm:infSuper}, that infinitely many polyhedra are superintegral, does not by itself imply \thmref{thm:infMany}.
\end{rmk}

Though 
every previously known integral packing was also superintegral, 
we discover for the first time that the latter is a strictly stronger condition.

\begin{lem}
There exist infnitely-many conformally-inequivalent crystallographic packings that are integral but not superintegral.
\end{lem}

\begin{rmk}\label{rmk:hexPyr}
Just one example of an integral but not superintegral 
polyhedron is
the hexagonal pyramid. See also \rmkref{rmk:nonArith}.
\end{rmk}

To prove \thmref{thm:infSuper}, we define certain 
operations on ``seed'' polyhedra which we call ``growths,''  including doubling 
the seed
along a vertex or a face, and observe that, while these generally wreak havoc on the resulting packings $\sP$, so $\sP(growth)$ and $\sP(seed)$ are usually  conformally inequivalent, 
 the superpackings $\widetilde\sP$ are essentially preserved, in fact
 $$ 
 \widetilde\sP(growth) \ \subset \ \widetilde\sP(seed).
 $$
 In particular, if a polyhedron is superintegral, then all of its growths are also superintegral, and hence integral! 
 This proves \thmref{thm:infSuper}, and hence 
\thmref{thm:infMany}.
\\

\noindent
{\bf \S(\ref{II}): Classifying Superintegral Crystallographic Packings}\

Towards the opposite general problem of classifying integral and superintegral crystallographic packings, we make
 two basic observations. The first, having nothing to do with integrality, shows that the entire theory of crystallographic packings is ``low''-dimensional.

\begin{thm}\label{thm:finDim}
Crystallographic packings can only exist in dimensions $n<996$.
\end{thm}

To prove this, we need the following
\begin{lem}\label{lem:lattice}
The supergroup $\widetilde\G$ of a crystallographic packing $\sP$ with symmetry group $\G$ is 
a {\em lattice}, that is, it acts on $\bH^{n+1}$ with {\it finite} covolume.
\end{lem}

We first sketch a proof of this lemma. Let $\G$ be a symmetry group for $\sP$; then it is assumed to be geometrically finite (recall that this means some uniform thickening of the convex core of $\G$ has finite volume). 
Since $\G$ is a reflection group, it has an essentially unique fundamental polyhedron $\sF:=\G\bk\bH^{n+1}$.
The domain of discontinuity $\gW$ of $\G$ (that is, the complement in $\dd\bH^{n+1}$ of its limit set $\gL_\G$) is the union of disjoint open geometric balls, since the limit set $\gL_\G$ is assumed to coincide with the set of limit points of $\sP$. 
The quotient $\gW/\G$ is then a disjoint union of finitely many open ends.
For each end, we 
develop the domain  under the $\G$-action  
and fill an open ball,
 the boundary of which is then an  (un-oriented) sphere in $\sP$. 
A geodesic hemisphere   
above such a ball is a frontier of the flare, cutting the walls it meets of 
$\sF$
either tangentially or at right angles (for otherwise the spheres in $\sP$ would overlap).
Hence when  we form the supergroup $\widetilde\G$ by adjoining  to $\G$ reflections through all the spheres in $\sP$,
we obtain a discrete action, and moreover the original domain of discontinuity $\gW$ has been entirely cut out, rendering $\widetilde\G$ a lattice. 
\\

Returning to \thmref{thm:finDim},
Vinberg \cite{Vinberg1981} and Prokhorov \cite{Prokhorov1986}  showed that hyperbolic reflection lattices  can only exist in dimensions $n<996$, and hence crystallographic packings are similarly bounded in dimension, proving the theorem. (The number $996$ is not expected to be sharp.) 
\\

Next we show that not only is the dimension bounded, but if we assume superintegrality, then (up to commensurability) there are only {\it finitely many}  Apollonian-like objects, period!

\begin{Def}
Two crystallographic packings are said to be {\it commensurate} if their supergroups are.
\end{Def}

\begin{thm}\label{thm:main}
There are only {\em finitely-many} commensurability classes 
 of superintegral crystallographic packings,  all of dimension $n<30$. 
\end{thm}


To prove this theorem, we show the following

\begin{thm}\label{conj:SuperPAC}
If $\sP$ is a superintegral crystallographic packing, then its supergroup $\widetilde\G$ is arithmetic\footnote{Recall that a real hyperbolic reflection
group is arithmetic if it is commensurate with the automorphism group of a hyperbolic quadratic form over the ring of integers of a totally real number field; see, e.g., \cite{Belolipetsky2016}.}!
\end{thm}

In fact, to conclude arithmeticity,  it is sufficient that the orbit under the supergroup $\widetilde\G$ of a single sphere $S\in\sP$ has all integer bends. 
Let us sketch a proof.
To a (positively-oriented) sphere $S$ of center $\bz=(z_1,\dots z_n)$ and radius $r$, we attach the ``inversive coordinates''
$$
\bv_S \ : = \ 
\left(\widehat b, b, b\bz\right).
$$
Here $b=1/r$ is the bend, and $\widehat b=1/\widehat r$ is the co-bend, that is, the reciprocal of the co-radius, the latter defined as the radius of the sphere after inversion through the unit sphere; see the discussion in, e.g., \cite{Kontorovich2017a, LagariasMallowsWilks2002}.
The vector $\bv_S$ lies on a one-sheeted hyperboloid $Q=-1$, where  $Q$ is the (universal) ``discriminant'' form, 
$$
Q=
\bp
 & \foh & \\
 \foh &  & \\
  &  & -I_{n-1}
\ep
.
$$
In these coordinates, 
\be\label{eq:GOQ}
\widetilde\G<O_Q(\R)
\ee 
is a right action by M\"obius transformations on the row vector $\bv_S$. Since $\widetilde\G$ is a lattice,
it is essentially (up to finite index components) Zariski dense in $O_Q$; hence  the orbit 
$
\cO=\bv_S\cdot\widetilde\G
$ 
of $S$  is essentially Zariski dense in the quadric $Q=-1$. There is then a choice of cluster $\cC_S\subset \cO$ of $n+2$ spheres whose matrix $\cV$ of inversive coordinates has (full) rank $n+2$. Make such a choice arbitrarily. This cluster $\cV$ has a Gram matrix of inversive products,
\be\label{eq:cGdef}
\cG \ := \ \cV\cdot Q\cdot \cV^\dag,
\ee
which is invertible (also has rank $n+2$). Let 
$$
\cF:=\cG^{-1}
$$
 be its inverse, which also induces a quadratic form having signature $(1,n+1)$. Then $\widetilde\G$ is conjugate to a  ``bends'' group, 
$$
\widetilde\cA \ := \ \cV\cdot\widetilde\G\cdot\cV^{-1} \ < \ O_\cF(\R)
,
$$ 
which now acts on the left on the (second) column vector of bends $\bb=\cV\cdot(0,1,0,\dots,0)^\dag$ in $\cV$; this 
vector $\bb$ lies on the  cone $\cF=0$, and $\widetilde\cA$ is a lattice in $O_\cF(\R)$.
Though {\it a priori} real valued, we claim that $\cF$ is in fact  {\it rational}.
Indeed, by assumption, the $\widetilde\cA$-orbit 
$$
\cB=\widetilde\cA\cdot\bb
$$
 lies in $\Z^{n+2}\cap \{\cF=0\}$,
and is Zariski dense in the cone. 
 But a quadratic form having a Zariski dense set of {\it integer} points $\cB$ on the cone $\cF=0$ is easily seen to be rational, as claimed. Next we observe that, since $\widetilde\cA$ is a {\it linear} action, it in fact preserves a full rank $\Z$-lattice $\gL$. But the group 
 $$
 O_\cF^\gL=\{g\in O_\cF(\R):g\gL=\gL\}
 $$ 
 is easily seen to be congruence, and contains $\widetilde\cA$. Hence $\widetilde\cA$ is arithmetic, as is its conjugate $\widetilde\G$. This proves \thmref{conj:SuperPAC}.
\\

Returning to  \thmref{thm:main}, this now follows from \thmref{conj:SuperPAC}, together with the amazing fact 
\cite{Vinberg1981, LongMaclachlanReid2006, Agol2006, ABSW2008, Nikulin2007}, 
that
there are only {\it finitely-many} 
commensurability classes of
arithmetic reflection groups, all having dimension $n< 30$. Hence the same holds for superintegral crystallographic packings
by \thmref{conj:SuperPAC}.
\\

It turns out that superintegrality is a necessary condition in \thmref{conj:SuperPAC}, and mere integrality is insufficient. Indeed, we discover for the first time the following
\begin{lem}\label{lem:23}
There exist infinitely many conformally inequivalent integral (but of course not superintegral) packings whose supergroups are non-arithmetic!
\end{lem}
\begin{rmk}\label{rmk:nonArith}
The supergroup of the hexagonal pyramid is {\it non-arithmetic};
see also \rmkref{rmk:hexPyr}.
\end{rmk}

\begin{rmk}
 Note also that there is no contradiction with \thmref{thm:infSuper} (and \thmref{thm:infMany}), as the 
 packings
 constructed there 
 fall
 into finitely many commensurability classes.
\end{rmk}

Given these finiteness results, the complete classification of superintegral crystallographic packings will then rely on understanding to what extent a converse of \thmref{conj:SuperPAC} may be true.

\begin{question}\label{q:conv}
Given an arithmetic reflection group, is it commensurate with the supergroup of a superintegral crystallographic packing?
\end{question}

We will say that an arithmetic group ``supports'' a packing if the answer to the above is YES.
We have investigated this question in some special cases and found the following positive results.

\begin{thm}\label{thm:17}
The answer to \qref{q:conv} is YES for all non-uniform lattices over $\Q$ in dimension $n=2$. Namely, every reflective (that is,  commensurate to a  reflection group)  Bianchi group supports a superintegral crystallographic packing.
\end{thm}

In higher dimensions, we are also able to show the following.

\begin{thm}\label{thm:highDim}
The answer is YES for certain lattices in all dimensions up to (at least) $n=13$; that is,
superintegral crystallographic packings exist in all these dimensions. 
\end{thm}

Before saying more about
 these theorems, let us point out that we suspect that the answer may be NO in general.
\begin{rmk}\label{rmk:coComp1}
At present, we do not know of a single superintegral (or even integral) packing whose supergroup is cocompact. In dimension $n=2$, the integral orthogonal groups preserving the form $x_1^2+x_2^2+x_3^2-dx_4^2$ are cocompact and reflective only when the coeffcient $d=7$ or $15$ \cite{McleodThesis}. We suspect, but do not know how to prove, that neither of these reflection groups support crystallographic packings. See \rmkref{rmk:coComp2}.
\end{rmk}
\begin{rmk}
 Taking, e.g.,  $\fo=\Z[\vf]$  the ring of the golden mean, we {\it can} construct $\fo$-superintegral packings (that is, with all bends in $\fo$), and having supergroup 
 the right-angled dodecahedron (which is arithmetic and co-compact).
It is an interesting problem to extend our theory to packings with bends in integer rings. (And more generally to complex hyperbolic space, $\SU(n,1)$, etc.)
\end{rmk}

\thmsref{thm:17} and \ref{thm:highDim} follow from our Structure Theorem: 

\begin{thm}[Structure Theorem for Crystallographic Packings]\label{thm:tildeC}
Let $\widetilde\cC$ be a  set of walls (that is, spheres), the reflections through which generate a hyperbolic lattice, and orient these walls
so that 
the fundamental domain
is the intersection of their exteriors.
Assume that $\widetilde\cC$ decomposes  into a cluster/cocluster pair:
\be\label{eq:tildeC}
\widetilde\cC=\cC\bigsqcup\widehat\cC
\ee
so that
\begin{itemize}
\item any pair of spheres in $\cC$ is either disjoint or tangent, and
\item
any sphere in $\cC$ is either disjoint, tangent, or orthogonal to any in $\widehat\cC$.
\end{itemize}
Let $\G:=\<\widehat\cC\>$ be the (thin) group generated by reflections through the cocluster.
Then the cluster orbit under this group, $\sP:=\G\cdot\cC$, is a crystallographic packing.

Conversely, every crystallographic packing arises in this way.
\end{thm}

The converse direction follows from our proof of \lemref{lem:lattice}, and the forward direction uses similar ideas.
Hence  answering  \qref{q:conv} for a given reflection lattice is equivalent to finding a decomposition as in \eqref{eq:tildeC}, or proving that one cannot exist.

\begin{rmk}\label{rmk:coComp2}
In the case of the cocompact forms in \rmkref{rmk:coComp1}, we are not yet able
after some effort
 to find a reflective subgroup (or prove it does not exist) with a suitable decomposition of the form \eqref{eq:tildeC}. 
 \end{rmk}






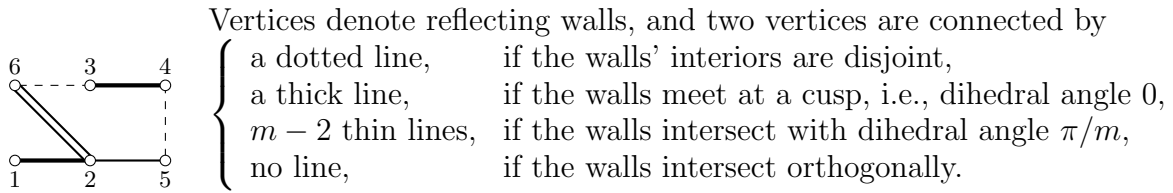
\begin{figure}

        \begin{subfigure}[t]{0.1\textwidth}
        \vskip.1in
 \begin{tikzpicture}
    \coordinate (one) at (-2, 0);
    \coordinate (two) at (-1, 0);
    \coordinate (three) at (-1, 1);
    \coordinate (four) at (0, 1);
    \coordinate (five) at (0, 0);
    \coordinate (six) at (-2, 1);

    \draw[ultra thick] (one) -- (two);
    \draw[thick] (two) -- (five);
    \draw[thick, double distance = 2pt] (two) -- (six);
    \draw[ultra thick] (three) -- (four);
    \draw[dashed] (three) -- (six);
    \draw[dashed] (four) -- (five);

{    \filldraw[fill=white] (one) circle (2pt) node[below] {\scriptsize 1};}
{    \filldraw[fill=white] (two) circle (2pt) node[below] {\scriptsize 2};}
{    \filldraw[fill=white] (three) circle (2pt) node[above] {\scriptsize 3};}
{    \filldraw[fill=white] (four) circle (2pt) node[above] {\scriptsize 4};}
{    \filldraw[fill=white] (five) circle (2pt) node[below] {\scriptsize 5};}
{    \filldraw[fill=white] (six) circle (2pt) node[above] {\scriptsize 6};}
  \end{tikzpicture}

        \end{subfigure}%
\hskip.5in
          \begin{subfigure}[t]{0.8\textwidth}
            Vertices denote reflecting walls, and two vertices are connected by

$
\left\{\begin{array}{ll}
        \text{a dotted line}, & \text{if the walls' interiors are disjoint,}\\
        \text{a thick line}, & \text{if the walls meet at a cusp, i.e., dihedral angle 0,}\\
        \text{$m-2$ thin lines}, & \text{if the walls intersect with dihedral angle $\pi/m$,}\\
        \text{no line}, & \text{if the walls intersect orthogonally.}\\
        \end{array}\right.
$
        \end{subfigure}

\caption{The Coxeter diagram for the reflective subgroup of the maximal discrete extension of the Bianchi group $\PSL_2(\Z[\sqrt{-6}])$.}
  \label{fig:Cox6}
\end{figure}

Returning to \thmref{thm:17}, our proof of this result relies on the complete classification by Belolipetsky-Mcleod \cite{BelolipetskyMcleod2013} of reflective Bianchi groups.
For example, 
the Bianchi group
$\PSL_2(\Z[\sqrt{-6}])$  
is commensurate to a
 maximal reflection group having the  Coxeter diagram illustrated in \figref{fig:Cox6}
%
(we follow Vinberg's convention for the labelling, as indicated there).

\begin{figure}
        \begin{subfigure}[t]{0.45\textwidth}
\includegraphics[width=.9\textwidth]{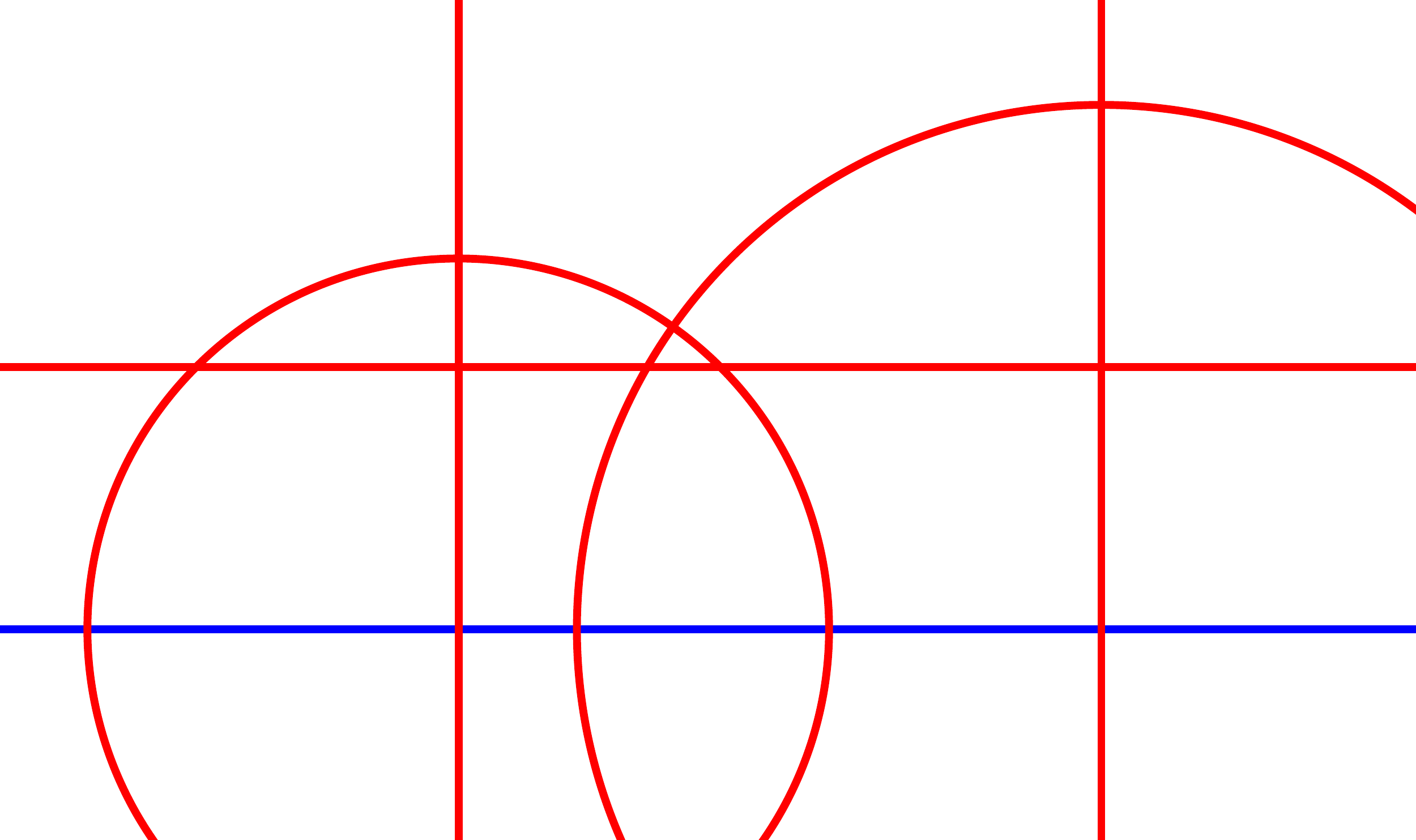}
\vskip-.4in
\hskip2.2in
$1$

\vskip-.65in
\hskip2.2in
$2$

\vskip-.8in
\hskip2.2in
$5$

\vskip-.2in
\hskip.65in
$4$

\vskip1in
\hskip.05in
$6$
\hskip1.75in
$3$

\vskip.1in
\caption{The reflecting walls in \figref{fig:Cox6}.}
  \label{fig:Pack6a}        
        \end{subfigure}%
          \begin{subfigure}[t]{0.45\textwidth}
\includegraphics[width=.9\textwidth]{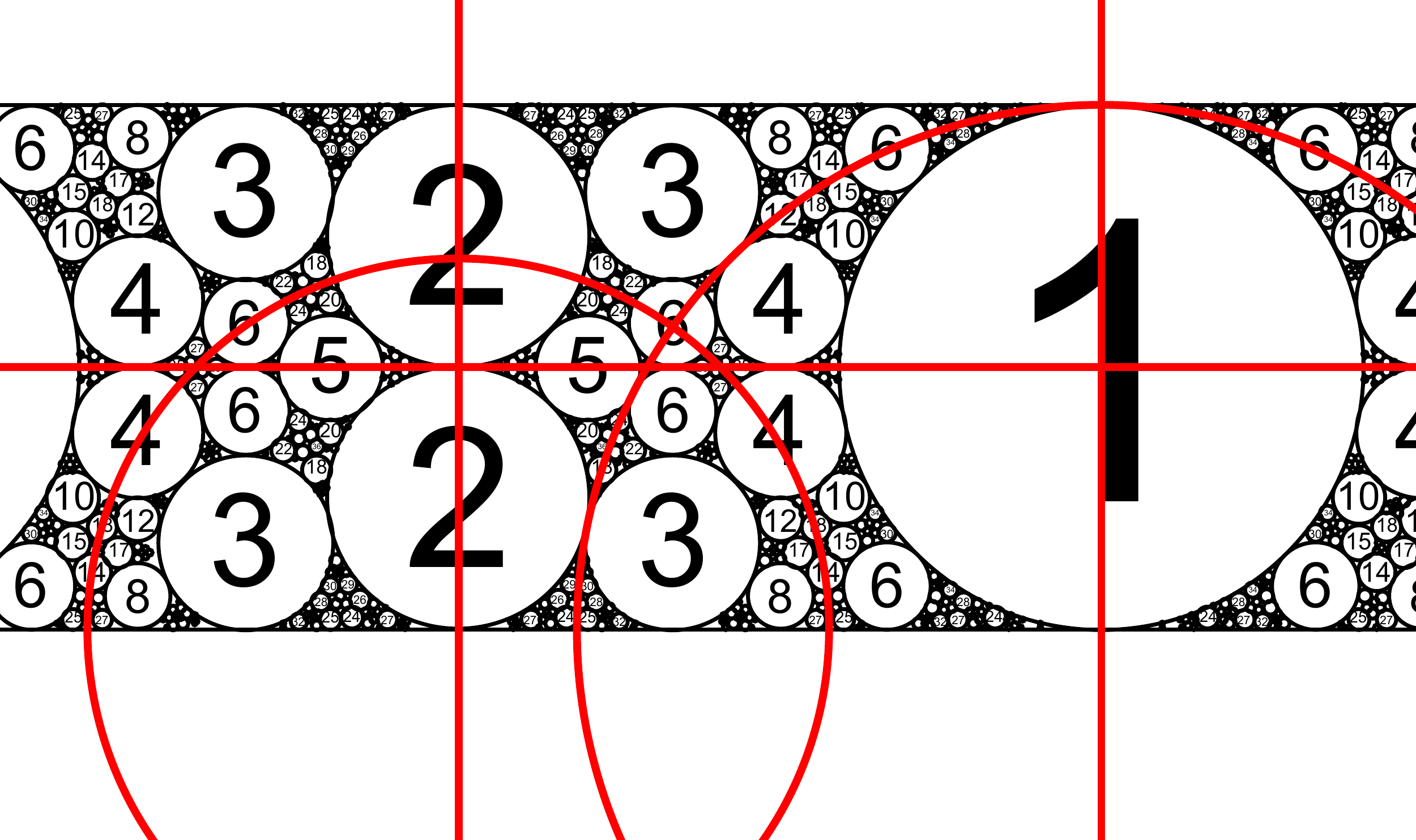}
\caption{The packing resulting from the orbit generated by $\widehat\cC$ on the cluster $\cC$ in the decomposition \eqref{eq:decompCubO}.}
  \label{fig:Pack6b}        
\end{subfigure}

\caption{}
  \label{fig:Pack6}
\end{figure}
One realization of the Coxeter diagram in \figref{fig:Cox6} is given by reflecting walls (circles) illustrated in \figref{fig:Pack6a}, with the same labeling. (The reader may check that the angles of intersection are as claimed in the Coxeter diagram.) 
The reader may now also verify that the decomposition of labeled walls as: 
\be\label{eq:decompCubO}
\cC=\{1\},\qquad \widehat\cC=\{2,3,\dots,6\},
\ee
satisfies the  
assumptions
of \thmref{thm:tildeC}, and
hence
 gives rise to a crystallographic packing by taking the orbit of $\cC$ under the group of reflections  through $\widehat\cC$.
The resulting packing is shown in \figref{fig:Pack6b},
which is the familiar cuboctahedral packing in disguise! (Compare with \figref{fig:cuboctPack}.)
\\

All but one of the other reflective Bianchi groups can be directly verified (by the \StructThm) to support superintegral packings; see Figures 1 and 2 in \cite{BelolipetskyMcleod2013}.
The only case in which the decomposition \eqref{eq:tildeC} is not straightforward from studying this Coxeter data is the Bianchi group on the Eisenstein  integers 
(that is, adjoining the cube root of unity).
It turns out in this case that the  Coxeter diagram in the literature has a minor mistake which can 
be traced to  an early paper of Shaiheev \cite{Shaiheev1990}; 
it has
 propagated in the literature ever since.
The issue  comes from the execution of Vinberg's algorithm for reflection subgroups which, for the Eisenstein integers, has extra stabilizers due to the larger 
group of units.
%
The 
true
diagram is

\begin{center}
  \begin{tikzpicture}
    \coordinate (one) at (-2, 0);
    \coordinate (two) at (-1, 0);
    \coordinate (three) at (0, 0);
    \coordinate (four) at (1, 0);

    \draw[thick] (four) -- (three);
    \draw[thick] (three) -- (two);
    \draw[thick, double distance = 2.5pt] (one) -- (two);
    \draw[thick, double distance = 0.3pt] (one) -- (two);

    \filldraw[fill=white] (one) circle (2pt) node[above] {\scriptsize 1};
    \filldraw[fill=white] (two) circle (2pt) node[above] {\scriptsize 2};
    \filldraw[fill=white] (three) circle (2pt) node[above] {\scriptsize 3};
    \filldraw[fill=white] (four) circle (2pt) node[above] {\scriptsize 4};
  \end{tikzpicture}
\end{center}

\noindent
and the Eisenstein Bianchi group  has a subgroup with Coxeter diagram

\begin{center}
  \begin{tikzpicture}
    \coordinate (one) at (-2, 0);
    \coordinate (five) at (-1, 0);
    \coordinate (three) at (2, 0);
    \coordinate (two) at (1, 0);
    \coordinate (four) at (0, 0);

    \draw[ultra thick] (one) -- (five);
    \draw[thick] (four) -- (five);
    \draw[ultra thick] (three) -- (two);
    \draw[thick, double distance = 2.5pt] (four) -- (two);
    \draw[thick, double distance = 0.3pt] (four) -- (two);

    \filldraw[fill=white] (one) circle (2pt) node[above] {\scriptsize 1};
    \filldraw[fill=white] (two) circle (2pt) node[above] {\scriptsize 4};
    \filldraw[fill=white] (three) circle (2pt) node[above] {\scriptsize 3};
    \filldraw[fill=white] (four) circle (2pt) node[above] {\scriptsize 5};
    \filldraw[fill=white] (five) circle (2pt) node[above] {\scriptsize 2};
  \end{tikzpicture}
\end{center}

\noindent
This last diagram
 supports a decomposition
 as in \eqref{eq:tildeC}  by taking either $\cC=\{1\}$ or $\cC=\{3\}$. We are thus finished sketching  the only non-immediate case of \thmref{thm:17}.

\begin{rmk}
In fact it turns out that {\it all} previously known integral circle packings (and many new ones!) arise in this way as limit sets of thin subgroups of reflective Bianchi groups.
\end{rmk}

To prove \thmref{thm:highDim},
we simply
inspect
 Vinberg's Coxeter diagrams \cite{Vinberg1972} for $n\le13$ for the  reflective subgroup of the integer orthogonal group preserving the form $-2x_0^2+x_1^2+\cdots+x_{n+1}^2$ 
 and apply the \StructThm.
\\


\noindent
{\bf \S\ Integral but Non-Superintegral Packings}\

Let us say more about what happens in \rmksref{rmk:hexPyr} and \ref{rmk:nonArith}. 
When $\Pi$ is, e.g., the hexagonal pyramid, its supergroup $\widetilde\G=\<\cC,\widehat\cC\>$ can be computed to have  Gram matrix (see \eqref{eq:cGdef})
\be\label{eq:Ghex}
\cG=
{
\left(
\begin{smallmatrix}
 -1 & 1 & 1 & 1 & 1 & 1 & 1 & 0 & 0 & 0 & 0 & 0 & 0 & \frac{2}{\sqrt{3}} \\
 1 & -1 & 1 & 5 & 7 & 5 & 1 & 0 & 2 \sqrt{3} & 4 \sqrt{3} & 4 \sqrt{3} & 2 \sqrt{3} & 0 & 0 \\
 1 & 1 & -1 & 1 & 5 & 7 & 5 & 0 & 0 & 2 \sqrt{3} & 4 \sqrt{3} & 4 \sqrt{3} & 2 \sqrt{3} & 0 \\
 1 & 5 & 1 & -1 & 1 & 5 & 7 & 2 \sqrt{3} & 0 & 0 & 2 \sqrt{3} & 4 \sqrt{3} & 4 \sqrt{3} & 0 \\
 1 & 7 & 5 & 1 & -1 & 1 & 5 & 4 \sqrt{3} & 2 \sqrt{3} & 0 & 0 & 2 \sqrt{3} & 4 \sqrt{3} & 0 \\
 1 & 5 & 7 & 5 & 1 & -1 & 1 & 4 \sqrt{3} & 4 \sqrt{3} & 2 \sqrt{3} & 0 & 0 & 2 \sqrt{3} & 0 \\
 1 & 1 & 5 & 7 & 5 & 1 & -1 & 2 \sqrt{3} & 4 \sqrt{3} & 4 \sqrt{3} & 2 \sqrt{3} & 0 & 0 & 0 \\
 0 & 0 & 0 & 2 \sqrt{3} & 4 \sqrt{3} & 4 \sqrt{3} & 2 \sqrt{3} & -1 & 1 & 5 & 7 & 5 & 1 & 1 \\
 0 & 2 \sqrt{3} & 0 & 0 & 2 \sqrt{3} & 4 \sqrt{3} & 4 \sqrt{3} & 1 & -1 & 1 & 5 & 7 & 5 & 1 \\
 0 & 4 \sqrt{3} & 2 \sqrt{3} & 0 & 0 & 2 \sqrt{3} & 4 \sqrt{3} & 5 & 1 & -1 & 1 & 5 & 7 & 1 \\
 0 & 4 \sqrt{3} & 4 \sqrt{3} & 2 \sqrt{3} & 0 & 0 & 2 \sqrt{3} & 7 & 5 & 1 & -1 & 1 & 5 & 1 \\
 0 & 2 \sqrt{3} & 4 \sqrt{3} & 4 \sqrt{3} & 2 \sqrt{3} & 0 & 0 & 5 & 7 & 5 & 1 & -1 & 1 & 1 \\
 0 & 0 & 2 \sqrt{3} & 4 \sqrt{3} & 4 \sqrt{3} & 2 \sqrt{3} & 0 & 1 & 5 & 7 & 5 & 1 & -1 & 1 \\
 \frac{2}{\sqrt{3}} & 0 & 0 & 0 & 0 & 0 & 0 & 1 & 1 & 1 & 1 & 1 & 1 & -1 \\
\end{smallmatrix}
\right)
}
.
\ee

Vinberg's Arithmeticity Criterion \cite{Vinberg1967} (see also \cite[Theorem 3.1]{VinbergShvartsman1993}) says in this context that $\widetilde\G$ is arithmetic if and only if cyclic products of $2\cG$ are always integers. 
This is almost the case for \eqref{eq:Ghex}, except for the entries $\frac2{\sqrt3}$ in the top right (and by symmetry, bottom left); hence $\widetilde\G$ is non-arithmetic (see \lemref{lem:23}). But it is nearly so; indeed, $\widetilde\G$, viewed as a subgroup of $O_Q(\R)$ (see \eqref{eq:GOQ}), can be conjugated to lie in $O_Q(\Z[\frac13])$ with unbounded denominators in its entries. The latter group is a perfectly nice $S$-arithmetic lattice in the product $O_Q(\R)\times O_Q(\Q_3)$, but $\widetilde\G$ is already a lattice on projection the first factor, $O_Q(\R)$. This too implies that $\widetilde\G$ is non-arithmetic, and in this sense is reminiscent of constructions of non-arithmetic groups by Deligne-Mostow \cite{DeligneMostow1986}. It is interesting to understand if all integral but non-superintegral packings arise this way. 
\\

\noindent
{\bf\S\  Local-Global Principles}\

We conclude with a discussion of whether  Local-Global Principles hold for bends of crystallographic circle ($n=2$) packings. (For higher dimensional sphere packings, this problem becomes easier; see, e.g., \cite{Kontorovich2017}.) As explained in \cite{Kontorovich2013} in the case of the classical Apollonian packing, the ``asymptotic'' local-global principle is proved in \cite{BourgainKontorovich2014b}.
This method was extended in the thesis of Zhang \cite{Zhang2015} to show the same statement for packings modeled on the octahedron. Most recently,
 Fuchs-Stage-Zhang showed that the method extends to the following context:
\begin{thm}[\cite{FuchsStangeZhang2017}]\label{thm:FSZ}
 Let $\sP$ be a packing with symmetry group $\G$ and let $C\in\sP$. Assume that  there is a circle $C'\in\sP$ tangent to $C$ so that the stabilizer of $C'$ in $\G$ is a \emph{congruence} (Fuchsian) group. Then the orbit $\G\cdot C$ satisfies an asymptotic local-global principle.
\end{thm} 

The assumption of the existence of such a companion circle $C'$ is a generalization of Sarnak's observation \cite{SarnakToLagarias} in the classical Apollonian case that such leads to certain shifted binary quadratic forms representing bends in the orbit. We show that this condition is both satisfied and not satisfied infinitely often!

\begin{thm}
The assumptions (and hence conclusions) of \thmref{thm:FSZ} are satisfied for infinitely many conformally-inequivalent superintegral crystallographic packings. 
The same statement holds with ``are satisfied'' replaced by ``are \emph{not} satisfied.''
\end{thm}

Thus even the asymptotic local-global problem remains 
open
 in this generality.

\subsection*{Acknowledegements}\

The authors benefitted
tremendously from numerous 
enlightening conversations about this work with
 Arthur Baragar, Misha Belolipetsky, Elena Fuchs, Jeremy Kahn, Jeff Lagarias, Alan Reid, Igor Rivin, Peter Sarnak, Kate Stange, Akshay Venkatesh,  Alex Wright, and Xin Zhang.

\bibliographystyle{alpha}

\bibliography{AKbibliog}

\end{document}